# Pinning Complex Networks by a Single Controller

Tianping Chen, Xiwei Liu, and Wenlian Lu

*Abstract*—In this paper, without assuming symmetry, irreducibility, or linearity of the couplings, we prove that a single controller can pin a coupled complex network to a homogenous solution. Sufficient conditions are presented to guarantee the convergence of the pinning process locally and globally. An effective approach to adapt the coupling strength is proposed. Several numerical simulations are given to verify our theoretical analysis.

## I. INTRODUCTION

MANY natural and man-made systems, such as neural systems, social systems, WWW, food webs, electrical power grids, etc., can all be described by graphs. In such a graph, every node represents an individual element of the system, while edges represent connections between nodes. Complex networks are such graphs with large size and complex topology. For decades, complex networks have been a focus on by scientists from various fields, for instance, sociology, biology, mathematics and physics, etc [1]–[3].

Linearly coupled ordinary differential equations (LCODEs) are a large class of dynamical systems with continuous time and state, as well as discrete space, which are widely used to describe coupling oscillators [4]. In general, the LCODEs can be described as follows:

$$\frac{dx_i(t)}{dt} = f(x_i(t),t) + c\sum_{j=1,j\neq i}^{m} a_{ij}[x_j(t) - x_i(t)],$$
$$i = 1,\ldots,m$$

where $x_i(t) = [x_i^1(t),\ldots,x_i^n(t)]^\top \in R^n$ is the state variable of the $i$th node, $t \in [0,+\infty)$ is the continuous time, $f : R^n \times [0,+\infty) \to R^n$ is continuous map, $a_{ij} \geq 0$ for $i,j = 1,\ldots,m$, $i \neq j$, denote the coupling coefficients, and $c$ is the coupling strength.

Letting $a_{ii} = -\sum_{j=1,j\neq i}^{m} a_{ij}$, the equations above can be rewritten as follows:

$$\frac{dx_i(t)}{dt} = f(x_i(t),t) + c\sum_{j=1}^{m} a_{ij}x_j(t), \qquad i = 1,\ldots,m \quad (1)$$

Manuscript received February 17, 2006; revised July 10, 2006. This work was supported by the National Science Foundation of China under Grant 60374018 and Grant 60574044.
T. Chen and X. Liu are with the Laboratory of Nonlinear Mathematics Science, Institute of Mathematics, Fudan University, Shanghai 200433, China (e-mail:tchen@fudan.edu.cn).
W. Lu is with the Laboratory of Nonlinear Mathematics Science, Institute of Mathematics, Fudan University, Shanghai 200433, China, and also with Max Planck Institute for Mathematics in the Sciences, Leipzig 04103, Germany.
Digital Object Identifier 10.1109/TCSI.2007.895383

where $A = (a_{ij}) \in R^{m \times m}$ is the coupling matrix with zero-sum rows and $a_{ij} \geq 0$, for $i \neq j$, which is determined by the topological structure of the LCODEs.

In past few decades, synchronization and control problems are being widely studied in complex networks. In [5]–[9], the local stability of the synchronization manifold was studied via the transverse stability to the synchronization manifold. The synchronizability based on the topology of the complex network was discussed in detail especially focusing on the complex networks with small-world and scale-free properties. In [10]–[12], a distance was defined from the collective spatial states of the coupled system to the synchronization manifold. Using this distance, methodologies were proposed to discuss global convergence for complete regular coupling configuration. In particular, in [13], the author pointed out that chaos synchronization can be obtained if and only if the topology of the network has a spanning tree.

Also, the problem of chaos control has been a research subject, which attracts increasing attention (see [14]–[16]for references). Recently, the object of chaos control has been transferred from single or several nodes to a dynamical networks especially complex networks (see [17], [18]). In particular, in [19]–[22], the authors studied pinning control problem on dynamical networks. Namely, controllers are only pinned on a very few fraction of nodes.

The problems of control and synchronization of coupled oscillators share similar background and have clear distinctiveness. In mathematical terms, (complete) synchronization can be described as follows: if a coupled system composed of $m$ sub-systems

$$\dot{x}_i = f_i(x_1,\ldots,x_m,t), \qquad i = 1,\ldots,m$$

satisfies $\lim_{t\to\infty}\|x_i(t) - x_j(t)\| = 0$, for all $i,j = 1,\ldots,m$, then the coupled system is said to be completely synchronized, for simplicity, synchronized. Instead, the pinning control problem is to synchronize all the states of the nodes in the dynamical network to a special solution $s(t)$ of the homogenous system $\dot{s}(t) = f(s(t))$. Namely, consider the following coupled system:

$$\dot{x}_i = f_i(x_1,\ldots,x_m,t) + \epsilon g_i(u_i), \qquad i \in J$$
$$\dot{x}_{i'} = f_{i'}(x_1,\ldots,x_m,t), \qquad i' \notin J$$

where $J$ is a subset of $\{1,\ldots,m\}$ denoting the controlled node set, $u_i$ denotes the control on the node $i \in J$, $g_i(\cdot)$ is some function realizing this control, and $\epsilon$ is the control strength. We define the pinning control performance by $\lim_{t\to\infty}\|x_i(t) - s(t)\| = 0$ holds for all $i = 1,\ldots,m$.

In [20], [21], the authors investigated pinning control for linearly coupled networks and found that one can pin the coupled networks by introducing fewer locally negative feedback controllers. They also compared two different pinning strategies: randomly pinning and selective pinning based on the connection degrees, and found out that the pinning strategy based on highest connection degree has better performance than totally randomly pinning. However, these studies were based on the intuitive knowledge of the network topologies, for example the connection degree. Instead, not based on rigorous mathematical analysis.

Therefore, it is natural to raise the following problem: **without any prior knowledge of the structure of the coupled network topology, can we pin the coupled network by introducing a single negative feedback controller (the simplest control)?**

In this paper, continuing with previous works, we study the pinning control problem of the coupled networks via a single controller. Without assuming symmetry, irreducibility, we discuss linearly or nonlinearly coupled networks. We provide sufficient conditions guaranteeing synchronizing a dynamical network to a homogenous solution with a single pinning controller. These criteria ensure that we can pin connected indirect or direct graphs with a spanning tree via a single controller. Moreover, we also propose an approach to adapt the coupling strength of the coupled network, which can significantly lessen the coupling strength.

This paper is organized as follows. In Section II, we present the main theoretical analysis for pinning linearly coupled complex networks via a single controller. In Section III, we discuss pinning nonlinearly coupled complex networks via a single controller. In Section IV, we provide an approach to adapt the coupling strength of the coupled network. Simulations are given in Section V to verify our theoretical results. We discuss and conclude the paper in Section VI and VII, respectively.

## II. Pin Linearly Coupled Complex Networks

Suppose that the linearly coupled network is

$$\frac{dx_i(t)}{dt} = f(x_i(t), t) + c\sum_{j=1}^{m} a_{ij} x_j(t) \quad i = 1, \ldots, m \quad (2)$$

where $x_i \in R^n$, $a_{ij} \geq 0$, $i \neq j$ and $\sum_{j=1}^{m} a_{ij} = 0$, for $i = 1, 2, \ldots, m$. $s(t)$ is a solution of the uncoupled system

$$\dot{s}(t) = f(s(t), t). \quad (3)$$

We prove that if $\varepsilon > 0$, the coupled network with a single controller shown in (4) at the bottom of the page, can pin the complex dynamical network (1) to $s(t)$, if $c$ is chosen suitably.

Denote $\delta x_i(t) = x_i(t) - s(t)$, then the system (1) can be rewritten as

$$\frac{d\delta x_i(t)}{dt} = f(x_i(t), t) - f(s(t), t) + c\sum_{j=1}^{m} a_{ij} \delta x_j(t),$$
$$i = 1, \ldots, m \quad (5)$$

and the network with a single controller (4) is written as

$$\frac{d\delta x_i(t)}{dt} = f(x_i(t), t) - f(s(t)) + c\sum_{j=1}^{m} \tilde{a}_{ij} \delta x_j(t),$$
$$i = 1, \ldots, m \quad (6)$$

where $\tilde{a}_{11} = a_{11} - \varepsilon$, $\varepsilon > 0$ and $\tilde{a}_{ij} = a_{ij}$ otherwise.

At first, we prove the following simple proposition.

*Proposition 1:* If $A = (a_{ij})_{i,j=1}^{m}$ is an irreducible matrix with Rank$(A) = m - 1$ and satisfying $a_{ij} = a_{ji} \geq 0$, if $i \neq j$, and $\sum_{j=1}^{m} a_{ij} = 0$, for $i = 1, 2, \ldots, m$. Then, all eigenvalues of the matrix

$$\tilde{A} = \begin{pmatrix} a_{11} - \varepsilon & a_{12} & \cdots & a_{1m} \\ a_{21} & a_{22} & \cdots & a_{2m} \\ \vdots & \vdots & \ddots & \vdots \\ a_{m1} & a_{m2} & \cdots & a_{mm} \end{pmatrix}$$

are negative.

*Proof:* Because $A$ is irreducible, there at least exists a positive element in the first column. Without loss of generality, we can assume $a_{21} > 0$. Let $\tilde{A}_1$ be the matrix obtained by excluding the first row and first column of $\tilde{A}$. Then, $\tilde{A}_1$ has the same structure with $\tilde{A}$.

Suppose that $\lambda$ is an eigenvalue of $\tilde{A}$, $v = [v_1, \ldots, v_m]^T$ is the corresponding eigenvector, and $|v_k| = \max_{j=1,\ldots,m} |v_j|$. It is clear that if $v$ is an eigenvector, then $-v$ is also an eigenvector. thus, without loss of generality, we can assume that $v_k > 0$ and $v_k = \max_{j=1,\ldots,m} |v_j|$.

If $k = 1$. Then

$$\sum_{j=1}^{m} \tilde{a}_{1j} v_j = -\varepsilon v_1 + \sum_{j=1}^{m} a_{1j} v_j$$
$$\leq -\varepsilon v_1 + \sum_{j=1}^{m} a_{1j} |v_j| < -\varepsilon v_1 < 0$$

which means $\lambda < 0$.

Instead, if $k > 1$. Then

$$\lambda v_k = \sum_{j=1}^{m} \tilde{a}_{kj} v_j \leq \tilde{a}_{kk} v_k + \sum_{j \neq k}^{m} \tilde{a}_{kj} |v_j| \leq 0$$

which means $\lambda \leq 0$.

$$\begin{cases} \frac{dx_1(t)}{dt} = f(x_1(t), t) + c\sum_{j=1}^{m} a_{1j} x_j(t) - c\varepsilon(x_1(t) - s(t)) \\ \frac{dx_i(t)}{dt} = f(x_i(t), t) + c\sum_{j=1}^{m} a_{ij} x_j(t), \quad i = 2, \ldots, m \end{cases} \quad (4)$$

If $\lambda = 0$. By the assumption that $A = (a_{ij})_{i,j=1}^m$ is an irreducible matrix with $\text{Rank}(A) = m - 1$. We conclude that $v = [v_k \ldots, v_k]^T$. However, this is impossible. For

$$\sum_{j=1}^m \tilde{a}_{1j} v_j = -\epsilon v_k < 0.$$

Therefore, $\lambda < 0$. The proposition is proved.

### A. Pin a Linearly Coupled Network With an Irreducible Symmetric Coupling Matrix

In this subsection, we investigate the linearly coupled networks, where the coupling matrix is symmetric.

With Proposition 1 given above, we prove two theorems. Theorem 1 addresses local synchronization. Theorem 2 addresses global synchronization.

Let $\delta x(t) = [\delta x_1(t), \ldots, \delta x_m(t)]$. Differentiating along $s(t)$ gives

$$\frac{d\delta x(t)}{dt} = D(f(s(t)), t)\delta x(t) + c\delta x(t)\tilde{A}^\top. \quad (7)$$

Let $\tilde{A}^T = WJW^\top$ be the eigenvalue decomposition of $\tilde{A}$, where $J = \text{diag}\{\lambda_1, \ldots, \lambda_m\}$, $0 > \lambda_1 > \cdots > \lambda_m$, and $\delta y(t) = \delta x(t)W$. Then we have

$$\frac{d\delta y_k(t)}{dt} = [Df(s(t),t) + c\lambda_k I]\delta y_k(t) \quad (8)$$

*Theorem 1:* Let $\mu_i(t)$, $i = 1, \ldots, m$, are the eigenvalues of the matrix $(1/2)(Df(s(t),t) + Df^T(s(t),t))$, $\mu(t) = \max_{i=1,\ldots,m} \mu_i(t)$. If $\mu(t) < -c\lambda_1 - \eta$ for all $t > 0$ and $\eta > 0$. Then, the coupled system with a controller (4) can be locally exponentially synchronized to $s(t)$.

*Proof:* It is easy to see that

$$\frac{1}{2}\frac{d\{\delta y_k^\top(t)\delta y_k(t)\}}{dt}$$
$$= \delta y_k^\top(t)[Df(s(t),t) + c\lambda_k I]\delta y_k(t)$$
$$= \delta y_k^\top(t)$$
$$\times \left[\frac{1}{2}(Df(s(t),t) + Df^T(s(t),t)) + c\lambda_k I\right]\delta y_k(t). \quad (9)$$

Under $\mu(t) < -c\lambda_1$ for all $t > 0$, we have

$$\frac{1}{2}\frac{d\{\delta y_k^\top(t)\delta y_k(t)\}}{dt} \le \delta y_k^\top(t)[\mu(t) - c\lambda_1]\delta y_k(t) \quad (10)$$

which means that $\delta y_k^\top(t)\delta y_k(t) = O(e^{-2\eta t})$. Theorem 1 is proved.

*Theorem 2:* Suppose $0 > \lambda_1 > \cdots > \lambda_m$ are the eigenvalues of $\tilde{A}$. If there are positive diagonal matrices $P = \text{diag}\{p_1, \ldots, p_n\}$, $\Delta = \text{diag}\{\Delta_1, \ldots, \Delta_n\}$ and a constant $\eta > 0$, such that

$$(x-y)^\top P(f(x,t) - \Delta x - f(y,t) + \Delta y)$$
$$\le -\eta(x-y)^\top(x-y) \quad (11)$$

and $\Delta_k + c\lambda_1 < 0$ for $k = 1, \ldots, n$. Then, the controlled system (4) is globally exponentially synchronized to $s(t)$.

*Proof:* Define a Lyapunov function as

$$V(t) = \frac{1}{2}\sum_{i=1}^m \delta x_i(t)^\top P\delta x_i(t).$$

Denote $\delta\tilde{x}^k(t) = [\delta x_1^k(t), \ldots, \delta x_m^k(t)]^\top$. Then, we have

$$\frac{dV(t)}{dt} = \sum_{i=1}^m \delta x_i(t)^\top Pf(x_i(t),t)\delta x_j(t)$$
$$- \sum_{i=1}^m \delta x_i(t)^\top Pf(s(t))\delta x_j(t)$$
$$+ \sum_{i=1}^m \delta x_i(t)^\top Pc \sum_{j=1}^m \tilde{a}_{ij}\delta x_j(t)$$
$$= \sum_{i=1}^m \delta x_i(t)^\top [P(f(x_i(t),t) - f(s(t))) - \Delta\delta x_i(t)]$$
$$+ \sum_{i=1}^m \delta x_i(t)^T P\left[c\sum_{j=1}^m \tilde{a}_{ij}\delta x_j(t) + \Delta\delta x_i(t)\right]$$
$$\le -\eta\sum_{i=1}^m \delta x_i(t)^\top \delta x_i(t)$$
$$+ \sum_{i=1}^m \delta x_i(t)^\top P\left[c\sum_{j=1}^m \tilde{a}_{ij}\delta x_j(t) + \Delta\delta x_i(t)\right]$$
$$= -\eta\sum_{i=1}^m \delta x_i(t)^\top \delta x_i(t)$$
$$+ \sum_{k=1}^n p_k\delta\tilde{x}^k(t)^\top(c\tilde{A} + \Delta_k I)\delta\tilde{x}^k(t).$$

Because $c\tilde{A} + \Delta_k I < 0$, we have

$$\frac{dV(t)}{dt} \le -\eta\sum_{i=1}^m \delta x_i(t)^\top \delta x_i(t) \le -2\eta\frac{V(t)}{\min_i p_i}.$$

Therefore

$$V(t) = O(e^{-2\eta t/\min_i p_i}).$$

Theorem 2 is proved completely.

*Remark 1:* It is clear that if $c$ is large enough, then the coupled network with a single controller can pin the complex network to a solution $s(t)$ of the uncoupled system.

*Remark 2:* Although the coupled network with a single controller can pin the complex network to a solution $s(t)$ of the uncoupled system. It does not mean that one must use one single controller to pin a coupled system. Theorem 2 also tells us that by adding any number of controllers can pin the coupled system. It is clear that the larger the number of the controllers is, the easier to pin a coupled system.

### B. Pin a Linearly Coupled Network With an Irreducible Asymmetric Coupling Matrix

In practice, indirect graphs are small part of the coupled networks. Most of the graphs are direct graphs. It means the coupling matrix is asymmetric. Therefore, we must investigate pining the complex networks, in which the coupling matrix is asymmetric. This is the issue we investigate in this section.

Consider the system

$$\frac{dx_1(t)}{dt} = f(x_1(t),t) + c\sum_{j=1}^m a_{1j}x_j(t) - c\varepsilon(x_1(t) - s(t))$$
$$\frac{dx_i(t)}{dt} = f(x_i(t),t) + c\sum_{j=1}^m a_{ij}x_j(t), \qquad i=2,\ldots,m \quad (12)$$

where the coupling matrix $A$ is asymmetric but satisfies zero row sum $\sum_{j=1}^m a_{ij} = 0, a_{ij} \geq 0$.

Let $[\xi_1,\ldots,\xi_m]^T$ be the left eigenvalue of the matrix $A$. It is well know that if $A$ is irreducible and $\text{Rank}(A) = m-1$, all $\xi_i > 0, i=1,\ldots,m$.

Define $\Xi = \text{diag}[\xi_1,\ldots,\xi_m]$. It is easy to verify that $\Xi A + A^T\Xi$ is a symmetric matrix and zero row sum. Therefore, by Proposition 1, the symmetric part of $\{\Xi\tilde{A}\}^s = (1/2)(\Xi\tilde{A} + \tilde{A}^T\Xi)$ is negative definite.

*Theorem 3:* Suppose that $A$ is irreducible and $\text{Rank}(A) = m-1$. $0 = \mu_1 > \mu_2 > \cdots > \mu_m$ are the eigenvalues of $\{\Xi\tilde{A}\}^s$. If there are positive diagonal matrices $P = \text{diag}\{p_1,\ldots,p_n\}$, $\Delta = \text{diag}\{\Delta_1,\ldots,\Delta_n\}$ and a constant $\eta > 0$, such that

$$(x-y)^\top P(f(x,t) - \Delta x - f(y,t) + \Delta y) \leq -\eta(x-y)^\top(x-y) \quad (13)$$

and $\Delta_k + c\mu_2 < 0$ for $k=1,\ldots,n$. Then, the controlled system (4) is globally exponentially synchronized to $s(t)$.

*Proof:* In this case, define a new Lyapunov function as

$$V(\delta x) = \frac{1}{2}\sum_{i=1}^m \xi_i \delta x_i^T P \delta x_i.$$

Differentiating $V(\delta x)$, we have

$$\frac{dV(\delta x)}{dt} = \sum_{i=1}^m \xi_i \delta x_i^T(t) P \frac{d\delta x_i(t)}{dt}$$
$$= \sum_{i=1}^m \xi_i \delta x_i^T(t) P[f(x_i(t),t) - f(s(t),t)]$$
$$+ c\sum_{i=1}^m \xi_i \delta x_i^T(t) P \sum_{j=1}^m \tilde{a}_{ij}\delta x_j(t)$$
$$= \sum_{i=1}^m \xi_i \delta x_i^T(t) P$$
$$\times \Bigg[(f(x_i(t),t) - \Delta x_i(t)) - (f(s(t),t)$$
$$- \Delta s(t)) + c\sum_{j=1}^m \tilde{a}_{ij}\delta x_j(t) + \Delta \delta x_i(t)\Bigg]$$

$$\leq -\eta\sum_{i=1}^m \xi_i \delta x_i^T(t)\delta x^i(t)$$
$$+ \sum_{i=1}^m \xi_i \delta x_i^T(t) P\Bigg[c\sum_{j=1}^m \tilde{a}_{ij}\delta x_j(t) + \Delta \delta x_i(t)\Bigg]$$
$$\leq -2\eta \frac{V(\delta x)}{\max_i p_i} + \sum_{j=1}^n p_j \delta \tilde{x}^{j\top}(t)\Xi(c\tilde{A} + \Delta_k I)\delta \tilde{x}^j(t)$$
$$= -2\eta \frac{V(\delta x)}{\max_i p_i} + \sum_{j=1}^n p_j \delta \tilde{x}^{j\top}(t)\left[\Xi(c\tilde{A} + \Delta_k I)\right]^s$$
$$\times \delta \tilde{x}^j(t) \leq -2\eta \frac{V(\delta x)}{\max_i p_i}.$$

For under the assumption $\Delta_j + c\mu_2 < 0$, we have $[\Xi(c\tilde{A} + \Delta_k I)]^s$ is negative definite for $k=1,\ldots,n$. Therefore

$$V(\delta x) = O(e^{-2\eta t/\max_i p_i}).$$

Theorem 3 is proved completely.

*C. Pin A Linearly Coupled Network With A Reducible Asymmetric Coupling Matrix*

In the following, we remove the assumption that $A$ is irreducible. In this case, we assume

$$A = \begin{bmatrix} A_{11} & 0 & \cdots & 0 \\ A_{21} & A_{22} & \cdots & 0 \\ \vdots & \ddots & \vdots & \vdots \\ A_{p1} & A_{p2} & \cdots & A_{pp} \end{bmatrix} \quad (14)$$

where $A_{qq} \in R^{m_q,m_q}$, $q=1,\ldots,p$, are irreducible or $A_{11}$ can be a zero matrix of one dimension. And, for each $q$, there exists $q > k$ such that $A_{qk} \neq 0$. It is equivalent to that the connecting graph has a spanning tree (see [13]).

It is easy to see that if we add a single controller $-\epsilon(x_1(t) - s(t))$ to the node $x_1(t)$. Then, by previous arguments, we conclude that the subsystem shown in (15) at the bottom of the page, pins $x_1(t),\ldots,x_{m_1}$ to $s(t)$.

Now, for the subsystem $x_{m_1+1},\ldots,x_{m_2}$, we have

$$\frac{d\delta x_i(t)}{dt} = f(x_i(t),t) - f(s(t),t)$$
$$+ c\sum_{j=1}^{m_2} a_{ij}\delta x_j(t)$$
$$= f(x_i(t),t) - f(s(t),t)$$
$$+ c\sum_{j=m_1+1}^{m_2} a_{ij}\delta x_j(t) + O(e^{-\eta t}).$$

---

$$\begin{cases} \frac{dx_1(t)}{dt} = f(x_1(t),t) + c\sum_{j=1}^{m_1} a_{1j}x_j(t) - c\varepsilon(x_1(t) - s(t)) \\ \frac{dx_i(t)}{dt} = f(x_i(t),t) + c\sum_{j=1}^{m_1} a_{ij}x_j(t), \qquad i=2,\ldots,m_1 \end{cases} \quad (15)$$

Because $A_{21} \neq 0$. Then, in $A_{22}$, there exists at least one row $i_2$, such that

$$a_{i_2 i_2} > \sum_{j=m_1+1}^{m_2} a_{i_2 j}. \qquad (16)$$

Therefore, all eigenvalues of $A_{22}$ are negative. By the similar arguments in the proof of Theorem 3, we can pin $x_i(t)$, $i = m_1 + 1, \ldots, m_2$, to $s(t)$.

By induction, we prove that if we add a controller to the master subsystem corresponding to the sub-matrix $A_{11}$, then we can ping the complex network to $s(t)$ even if the coupling matrix in reducible.

*Remark 3:* An important issue is the estimation of the upper-bound for the coupling strength. Because, in practice, it can not be realized if the coupling strength is too large. It also depends on the three factors in previous remark. In this paper, we will not discuss this issue. In Section III, we propose an adaptive approach to find the coupling strength, which is significantly smaller than the theoretical value.

## III. PIN A NONLINEARLY COUPLED NETWORK WITH A SYMMETRIC COUPLING MATRIX

In this section, we discuss how to pin a complex network with nonlinear coupling.

In many practical problems, it often happens that the states $x_i(t)$ can not be observed directly. Instead, we can only observe data $g(x_i(t)) = [g(x_i^1(t)), \ldots, g(x_i^n(t))]^\top$, $i = 1, \ldots, m$. We need to synchronize the uncoupled system by these data $g(x^i(t))$. It means that the synchronization scheme is nonlinear. Hence, investigation of synchronization for nonlinear coupled dynamical networks is an necessary step in both theoretical research and applications. Moreover, to describe $x_i(t)$ properly, every function $g_i(\cdot)$ should be monotone.

In this case, the coupled system (4) takes (17), shown at the bottom of the page, where $g(x_i(t)) = [g(x_i^1(t)), \ldots, g(x_i^n(t))]^\top$ and $g(x)$ is a monotone increasing function.

In the following, we will prove that the coupled complex network with a single controller (17) can be pined to a specified solution $\dot{s}(t) = f(s(t),t)$, too. In particular, we prove

*Theorem 4:* Suppose $0 > \lambda_1 > \cdots > \lambda_m$ are the eigenvalues of $\tilde{A}$, $(g(u) - g(v)/u - v) \geq \underline{\alpha} > 0$. If there are positive diagonal matrices $P = \text{diag}\{p_1, \ldots, p_n\}$, $\Delta = \text{diag}\{\Delta_1, \ldots, \Delta_n\}$ and a constant $\eta > 0$, such that

$$(x-y)^\top P(f(x,t) - \Delta x - f(y,t) + \Delta y) \leq -\eta(x-y)^\top(x-y) \qquad (18)$$

and $\Delta_k + \underline{\alpha} c \lambda_1 < 0$ for $k = 1, \ldots, n$. Then, the controlled system (17) is globally exponentially synchronized to $s(t)$.

*Proof:* Along with $g(x_i) = [g(x_i^1), \ldots, g(x_i^n)]^\top$, $i = 1, \ldots, m$, we denote $\tilde{g}(x^k) = [g(x_1^k), \ldots, g(x_m^k)]^\top$, $k = 1, \ldots, n$.

We use the same Lyapunov function

$$V(t) = \frac{1}{2} \sum_{i=1}^{m} \delta x_i(t)^\top P \delta x_i(t).$$

In this case, we have

$$\begin{aligned}\frac{dV(t)}{dt} &= \sum_{i=1}^{m} \delta x_i^\top(t) P \frac{d\delta x_i(t)}{dt} \\ &= \sum_{i=1}^{m} \delta x_i^\top(t) P \\ &\quad \times \left[f(x^i(t),t) - f(s(t),t) + c\sum_{j=1}^{m} a_{ij}g(x_j(t))\right] \\ &\quad - c\varepsilon[x_1(t) - s(t)]^\top P[g(x_1(t)) - g(s(t))] \\ &\leq -\eta \sum_{i=1}^{m} \delta x_i^\top(t)\delta x^i(t) + \sum_{i=1}^{m} \delta x_i^\top(t) P \\ &\quad \times \left[\Delta \delta x_i(t) + c\sum_{j=1}^{m} a_{ij}g(x_j(t))\right] \\ &\quad - c\varepsilon[x_1(t) - s(t)]^\top P \underline{\alpha}[x_1(t) - s(t)].\end{aligned}$$

By the property of the matrix $A$, it is easy to verify that for $u = [u_1, \ldots, u_m]^\top$, $v = [v_1, \ldots, v_m]^\top$,

$$u^\top A v = \sum_{i,j=1}^{m} u_i a_{ij} v_j = \sum_{j>i} a_{ij}(u_i - u_j)(v_i - v_j) \qquad (19)$$

combining with $(g(u) - g(v)/u - v) \geq \underline{\alpha} > 0$, we have

$$\begin{aligned}\sum_{i=1}^{m} &\delta x_i^\top(t) P \sum_{j=1}^{m} a_{ij}g(x_j(t)) \\ &= \sum_{k=1}^{n} p_k \delta \tilde{x}^k(t)^\top A \delta \tilde{g}(x^k(t)) \\ &= \sum_{k=1}^{n} p_k \sum_{j>i} a_{ij}(x_i^k(t) - x_j^k(t))(g(x_i^k(t)) - g(x_j^k(t))) \\ &\leq -\underline{\alpha} \sum_{k=1}^{n} p_k \sum_{j>i} a_{ij}(x_i^k(t) - x_j^k(t))(x_i^k(t) - x_j^k(t)) \\ &= -\underline{\alpha} \sum_{k=1}^{n} p_k \delta \tilde{x}^k(t)^\top A \delta \tilde{x}^k(t).\end{aligned}$$

$$\begin{cases} \frac{dx_1(t)}{dt} = f(x_1(t),t) + c\sum_{j=1}^{m} a_{1j}g(x_j(t)) - c\varepsilon(g(x_1(t)) - g(s(t))) \\ \frac{dx_i(t)}{dt} = f(x_i(t),t) + c\sum_{j=1}^{m} a_{ij}g(x_j(t)), \qquad i = 2,\ldots,m \end{cases} \qquad (17)$$

Therefore,

$$\frac{dV(t)}{dt} \leq -\eta \sum_{i=1}^{m} \delta x_i^\top(t)\delta x^i(t)$$
$$- \sum_{k=1}^{n} p_k \delta\tilde{x}^k(t)^\top [\Delta_k I_m + c\underline{\alpha}\tilde{A}]\delta\tilde{x}^k(t).$$

Because $\Delta_k + \underline{\alpha} c\lambda_1 < 0$, we have

$$\frac{dV(t)}{dt} \leq -2\eta \frac{V(t)}{\max_i p_i}$$
$$V(t) = O(e^{-2\eta t/\max_i p_i}).$$

Theorem 4 is proved.

## IV. ADAPTIVE ADJUSTMENT OF THE COUPLING STRENGTH

In the previous section, we proved that we can always pin a coupled complex network by adding a single controller if the coupling strength is large enough. However, in practice, it is not allowed that the coupling strength is arbitrarily large. For synchronization, it was pointed out in [11] that theoretical value of the coupling strength is much larger than needed in practice. Therefore, the following question was arisen in [11]: *Can we find the sharp bound $c_{\min}$* Similarly, in pinning process, it is also important to make the coupling strength as small as possible. It is clear that theoretical value of strength given in previous theorems are based on the condition (11)

$$(x-y)^\top P(f(x,t) - \Delta x - f(y,t) + \Delta y)$$
$$\leq -\eta (x-y)^\top (x-y) \quad (20)$$

which is too strong. In fact, for many chaotic oscillators, with many $x$, $y$, we even have

$$(x-y)^\top P(f(x,t) - f(y,t)) \leq -\eta(x-y)^\top(x-y). \quad (21)$$

Therefore, it is possible to lessen coupling strength dramatically.

In case we don't know the structure of the coupled system, a simple approach is to adapt the coupling strength. For this purpose, we prove

*Theorem 5:* Under the assumptions of Theorem 3, the coupled system as shown in (22), at the bottom of the page, where $c(0) \geq 0$ and $\alpha > 0$, can synchronize to the given trajectory $s(t)$.

*Proof:* Pick a constant $\alpha > 0$. Define a Lyapunov function

$$V(t) = \frac{1}{2}\sum_{i=1}^{m} \xi_i \delta x_i^T(t) P \delta x_i^T(t) + \frac{\beta}{\alpha}(c - c(t))^2; \quad (23)$$

where constants $c$ and $\beta$ will be decided later.

Differentiating it, we have

$$\frac{dV(t)}{dt} = \sum_{i=1}^{m} \xi_i \delta x_i^T(t) P$$
$$\times \left[ f(x_i(t),t) - f(s(t),t) + c(t)\sum_{j=1}^{m}\tilde{a}_{ij}\delta x_j(t) \right]$$
$$- \beta(c - c(t))\sum_{i=1}^{m} \delta x_i^T(t) P \delta x_i(t)$$
$$= \sum_{i=1}^{m} \xi_i \delta x_i^T(t) P \left[ f(x_i(t),t) - f(s(t),t) - \Delta \delta x_i(t) \right]$$
$$+ \sum_{i=1}^{m} \xi_i \delta x_i^T(t) P \Delta \delta x_i(t)$$
$$+ c(t)\sum_{i=1}^{m} \xi_i \delta x_i^T(t) P \sum_{j=1}^{m} \tilde{a}_{ij}\delta x_j(t)$$
$$- \beta(c - c(t))\sum_{i=1}^{m} \delta x_i^T(t) P \delta x_i(t)$$
$$\leq -\eta \sum_{i=1}^{m} \xi_i \delta x_i(t)^T \delta x_i(t) + \sum_{i=1}^{m} \xi_i \delta x_i^T(t) P \Delta \delta x_i(t)$$
$$+ c(t)\sum_{k=1}^{n} p_k \delta\tilde{x}_k^T(t) \Xi \tilde{A} \delta\tilde{x}_k(t)$$
$$- \beta c \sum_{i=1}^{m} \delta x_i^T(t) P \delta x_i(t) + \beta c(t)\sum_{i=1}^{m} \delta x_i^T(t) P \delta x_i(t)$$
$$\leq -\eta \sum_{i=1}^{m} \xi_i \delta x_i(t)^T \delta x_i(t)$$
$$+ \max_{l} \xi_l \sum_{k=1}^{n} p_k \Delta_k \delta\tilde{x}_k^T(t)\delta\tilde{x}_k(t)$$
$$+ c(t)\sum_{k=1}^{n} p_k \delta\tilde{x}_k^T(t)\Xi\tilde{A}\delta\tilde{x}_k(t)$$
$$- \beta c \sum_{k=1}^{n} p_k \delta\tilde{x}_k^T(t)\delta\tilde{x}_k(t)$$
$$+ \beta c(t)\sum_{k=1}^{n} p_k \delta\tilde{x}_k^T(t)\delta\tilde{x}_k(t)$$
$$= -\eta \sum_{i=1}^{m} \xi_i \delta x_i(t)^T \delta x_i(t)$$
$$+ \sum_{k=1}^{n} p_k (\max_l \xi_l \Delta_k - \beta c)\delta\tilde{x}_k^T(t)\delta\tilde{x}_k(t)$$
$$+ c(t)\sum_{k=1}^{n} p_k \delta\tilde{x}_k^T(t)(\Xi\tilde{A} + \beta I)\delta\tilde{x}_k(t).$$

$$\begin{cases} \frac{dx_1(t)}{dt} = f(x_1(t),t) + c(t)\sum_{j=1}^{m} a_{1j}x_j(t) - c(t)\varepsilon(x_1(t) - s(t)) \\ \frac{dx_i(t)}{dt} = f(x_i(t),t) + c(t)\sum_{j=1}^{m} a_{ij}x_j(t), \quad i=2,\ldots,m \quad \dot{c}(t) = \frac{\alpha}{2}\sum_{i=1}^{m}\delta x_i^T(t)P\delta x_i(t) \end{cases} \quad (22)$$

It is clear that $c(t) > 0$ for for all $t > 0$. Therefore, by picking constants $c$ and $\beta$ such that $\{\Xi\tilde{A} + \beta I\}^s < 0$ and $\max_{l,k} \xi_l \Delta_k - \beta c < 0$, then we have

$$\frac{dV(t)}{dt} < 0.$$

It is obvious that the set $H_1 = \{\delta x_i(t) = 0, i = 1, \ldots, m\}$ contains the largest invariant set contained in $H_2 = \{(dV/dt) = 0\}$. By the invariant principle of functional differential equations, we conclude $x_i(t) \to s(t)$ and from (22), when $x_i(t) \to s(t)$, then $\dot{c}(t) \to 0$, i.e., $c(t) \to c_0$, where $c_0$ is a positive constant. The proof is completed.

## V. NUMERICAL EXPERIMENTS

In this section, we give numerical examples to verify the theorems given in previous sections.

We consider chaotic Chen's oscillators, Lorenz oscillators, Rössler oscillators and Chua's circuits as the examples of the uncoupled network.

A single Lorenz oscillator is described in the dimensionless form by

$$\begin{cases} \frac{dx_1}{dt} = \beta[x_2 - x_1] \\ \frac{dx_2}{dt} = \alpha x_1 - x_2 - x_1 x_3 \\ \frac{dx_3}{dt} = x_1 x_2 - b x_3 \end{cases} \quad (24)$$

where $\beta = 10$, $\alpha = 28$, $b = (8/3)$

A single Chen's oscillator is described in the dimensionless form by

$$\begin{cases} \frac{dx_1}{dt} = a[x_2 - x_1] \\ \frac{dx_2}{dt} = (c - a)x_1 - x_2 x_3 + c x_2 \\ \frac{dx_3}{dt} = x_1 x_2 - b x_3 \end{cases} \quad (25)$$

where $a = 35$, $b = 3$, $c = 28$.

A single Rössler oscillator is described in the dimensionless form by

$$\begin{cases} \frac{dx_1}{dt} = -[x_2 - x_3] \\ \frac{dx_2}{dt} = x_1 + 0.2 x_2 \\ \frac{dx_3}{dt} = 0.2 + x_3(x_1 - \mu) \end{cases} \quad (26)$$

where $\mu = 5.7$.

A single Chua's circuit is described in the dimensionless form by

$$\begin{cases} \frac{dx_1}{dt} = k[x_2 - h(x_1)] \\ \frac{dx_2}{dt} = x_1 - x_2 + x_3 \\ \frac{dx_3}{dt} = -l x_2 \end{cases} \quad (27)$$

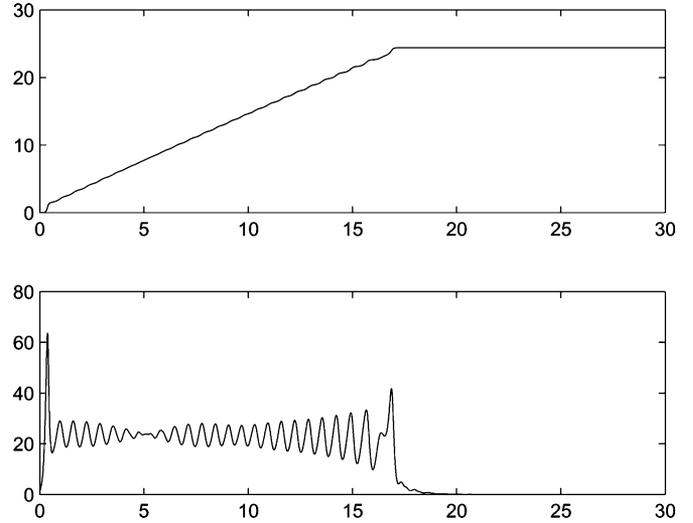

Fig. 1. Pin 500 Lorenz oscillators coupled by a small-world network with one controller.

where $h(x) = (2/7)x - (3/14)[|x+1| - |x-1|]$, $k = 9$ and $l = 14(2/7)$.

### A. Simulation 1

In this simulation, we consider a small-world network with 500 nodes. The connection weights are chosen randomly in $[0, 1]$ with uniform distribution.

We use the coupled network shown in (28) at the bottom of the page, to pin all nodes to a specified trajectory of the corresponding chaotic oscillators, where $f(\cdot)$ is chaotic Lorenz oscillator, Chen's chaotic oscillator, Rössler chaotic oscillator, and Chua's circuit, respectively. The index $i_0$ is chosen so that $\sum_{j \neq i_0} |a_{ji_0}| = \max_i \{\sum_{j \neq i} |a_{ji}|\}$. The control gain $\varepsilon$ is 100. The quantity $E(t) = \sqrt{\left(\sum_{i=1}^{500} |x_i(t) - s(t)|^2/500\right)}$ is used to measure the quality of the pinning process. In all the following figures, the upper one indicates the evolution of the coupling strength $c(t)$. The lower one indicates the evolution of $E(t)$.

Fig. 1 indicates how $E(t)$ and $c(t)$ evolute in pinning 500 Lorenz oscillators coupled by a small-world network with final coupling strength $C = 24.403998$.

Fig. 2 indicates how $E(t)$ and $c(t)$ evolute in pinning 500 coupled Chen's oscillators coupled by a small-world network with final coupling strength $C = 39.632460$.

Fig. 3 indicates how $E(t)$ and $c(t)$ evolute in pinning 500 Rössler oscillators coupled by a small-world network with final coupling strength $C = 6.180105$.

$$\begin{cases} \frac{dx_{i_0}(t)}{dt} = f(x_{i_0}(t), t) + c(t) \sum_{j=1}^{500} a_{i_0 j} x_j(t) - c(t)\varepsilon(x_{i_0}(t) - s(t)) \\ \frac{dx_i(t)}{dt} = f(x_i(t), t) + c(t) \sum_{j=1}^{500} a_{ij} x_j(t), \quad i = 2, \ldots, m, i \neq i_0 \quad \dot{c}(t) = \frac{\alpha}{2} \sum_{i=1}^{m} \delta x_i^T(t) P \delta x_i(t) \end{cases} \quad (28)$$

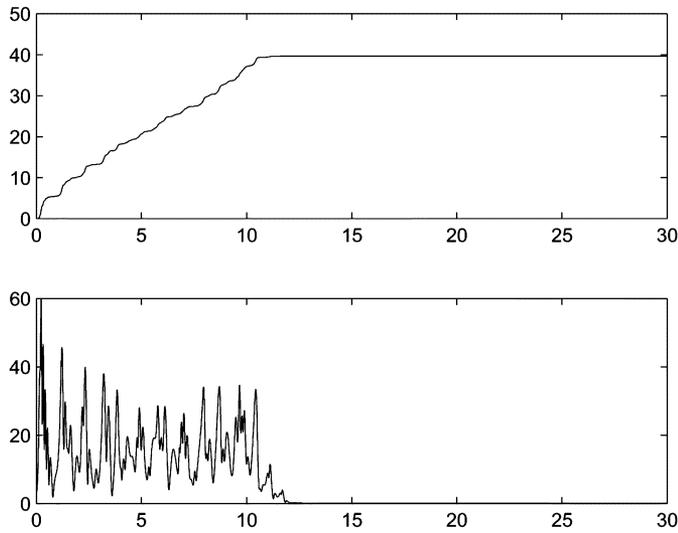

Fig. 2. Pin 500 Chen's oscillators coupled by a small-world network with one controller.

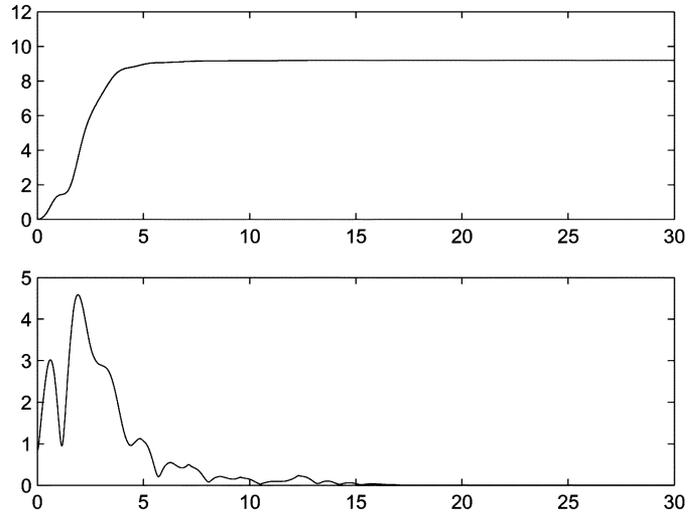

Fig. 4. Pin 500 Chua's circuits coupled by a small-world network with one controller.

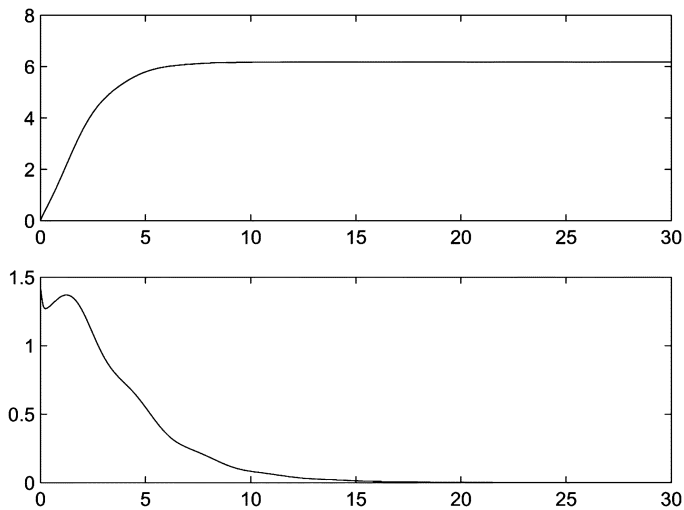

Fig. 3. Pin 500 Rössler oscillators coupled by a small-world network with one controller.

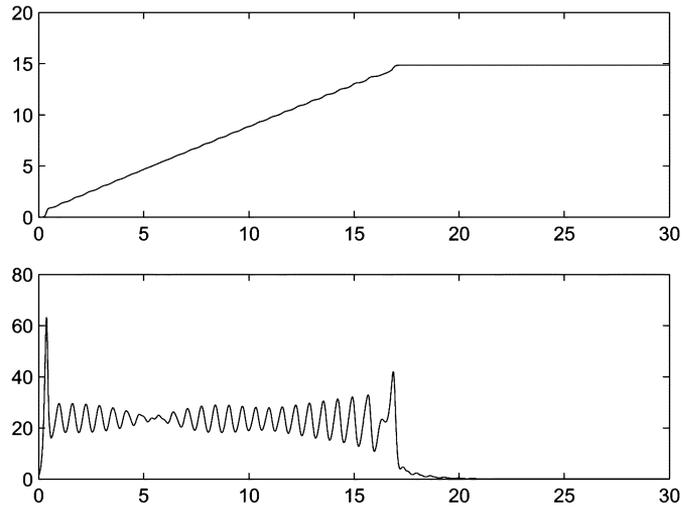

Fig. 5. Pin 500 Lorenz oscillators coupled by a randomly chosen network with one controller.

Fig. 4 indicates how $E(t)$ and $c(t)$ evolute in pinning 500 Chua's circuits coupled by a small-world network with final coupling strength $C = 9.185478$.

*Remark 4:* In [20], the authors pinned 50 Chen's oscillators coupled by a coupled small-world network by adding two controllers with control gain $\epsilon = 1000$. Here, we pin 500 oscillators by one controller with control gain $\epsilon = 100$. The coupling strength is adapted.

*Remark 5:* Placing the controller node may influence the coupling strength. However, in our many other simulations, this influence is not very heavy.

### B. Simulation 2

In this section, the coupling matrix has only 20 percent non-zero entries, which are chosen randomly in $[0,1]$ with uniform distribution. As the chaotic oscillators, we also use chaotic Lorenz oscillator, Chen's chaotic oscillator, Rössler chaotic oscillator and Chua's circuit, respectively.

Fig. 5 indicates how $E(t)$ and $c(t)$ evolute in pinning 500 Lorenz oscillators coupled by a randomly chosen network with final coupling strength $C = 14.853546$.

Fig. 6 indicates how $E(t)$ and $c(t)$ evolute in pinning 500 Chen's oscillators coupled by a randomly chosen network with final coupling strength $C = 37.831439$.

Fig. 7 indicates how $E(t)$ and $c(t)$ evolute in pinning 500 Rossler oscillators coupled by a randomly chosen network with final coupling strength $C = 5.327969$.

Fig. 8 indicates how $E(t)$ and $c(t)$ evolute in pinning 500 Chua's circuits oscillators coupled by a randomly chosen network with final coupling strength $C = 6.149759$.

At last, we give two simulation experiments. One is to pin 1000 Rössler oscillators coupled by a randomly chosen network. The other one is pin 1000 Chua's circuits coupled by a randomly chosen network.

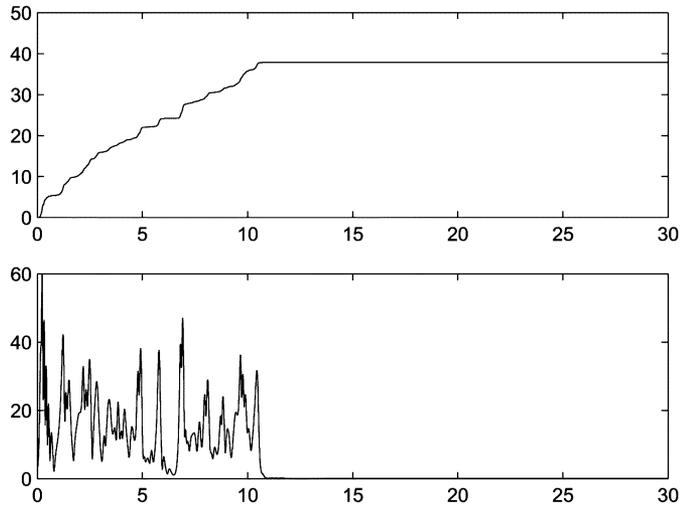

Fig. 6. Pin 500 Chen's oscillators coupled by a randomly chosen network with one controller.

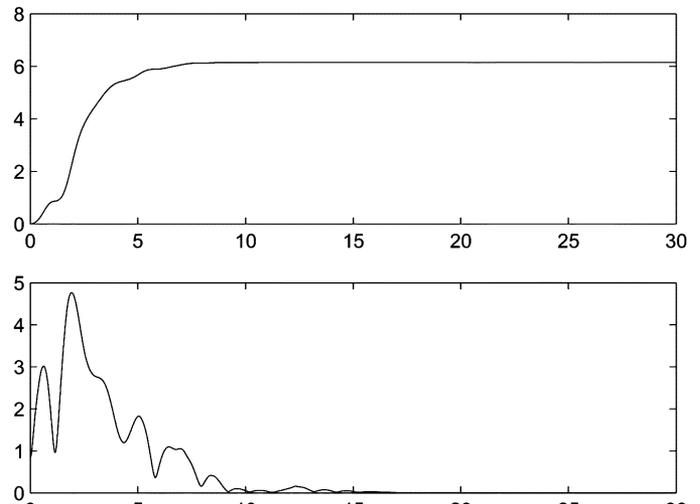

Fig. 8. Pin 500 Chua's circuits coupled by a randomly chosen network with one controller.

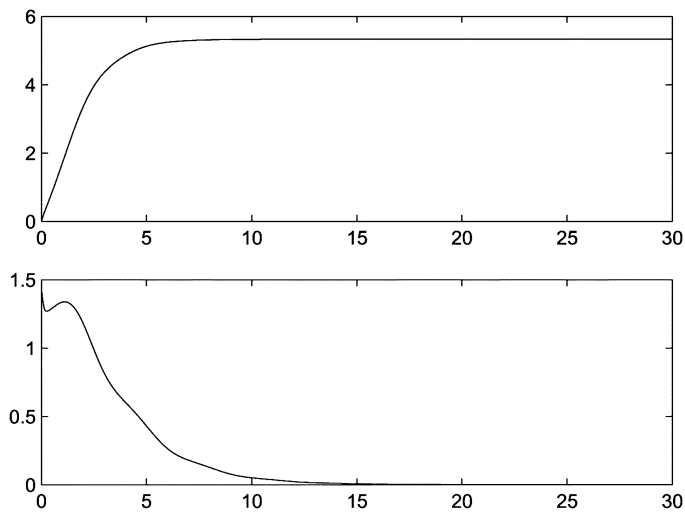

Fig. 7. Pin 500 Rössler oscillators coupled by a randomly chosen network with one controller.

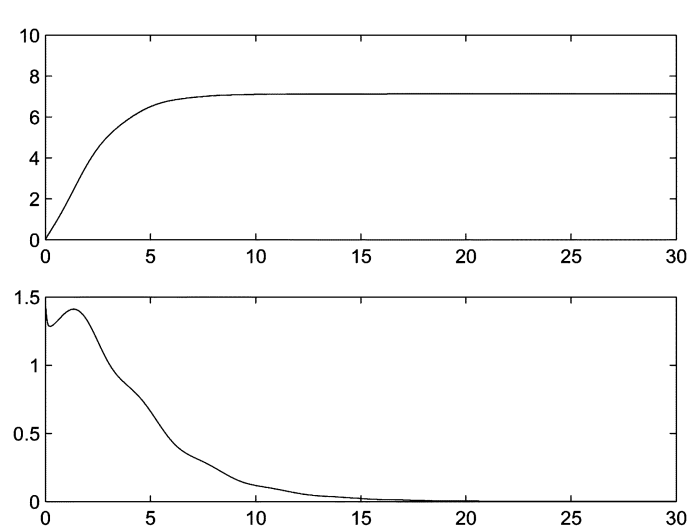

Fig. 9. Pin 1000 Rössler oscillators coupled by a randomly chosen network with one controller.

Fig. 9 indicates how $E(t)$ and $c(t)$ evolute in pinning 1000 Rössler oscillators coupled by a randomly chosen network with final coupling strength 7.124347.

Fig. 10 indicates how $E(t)$ and $c(t)$ evolute in pinning 1000 Chua's circuits coupled by a randomly chosen network with final coupling strength $C = 8.893303$.

*Remark 6:* It is very interesting that one can pin 1000 chaotic oscillators coupled by adding a single controller with quite small coupling strength.

*Remark 7:* It is thought that the coupling strength will increase as long as the number of the nodes coupled. However, in many simulations, the coupling strength changes not so dramatically when the number of the nodes increases. Instead, the accuracy of the algorithm to solve ordinary differential equations plays a key role. The more the number of nodes, the more accuracy is required.

*Remark 8:* The location of the pinned node does not influence the performance heavily as shown in Fig. 9. The node to add controller is chosen randomly.

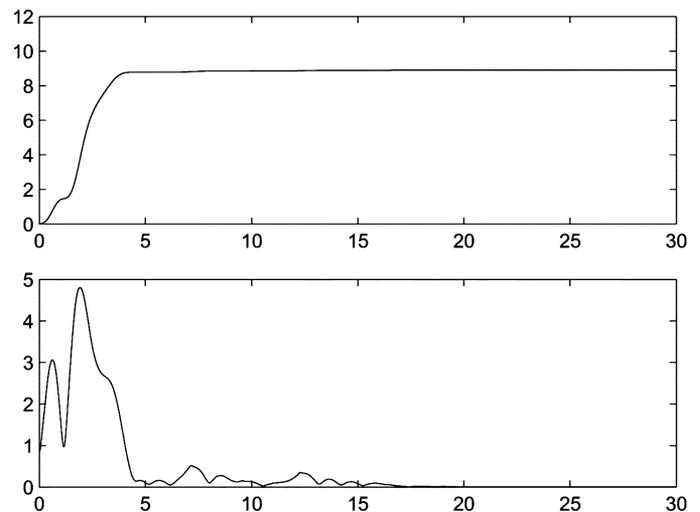

Fig. 10. Pin 1000 Chua's circuits coupled by a randomly chosen network with one controller.

## VI. Discussion

In this paper, we prove that one can pin a coupled system to a specified trajectory by inducing a single controller. This problem can also be rewritten in the synchronization framework as follows:

$$\dot{x}_i(t) = f(x_i(t),t) + c(t)\sum_{j=0}^{m} \tilde{a}_{ij} x_j(t), \qquad i=,0,1,\ldots,m \quad (29)$$

where $x_0(t) = s(t)$, $\tilde{a}_{0j} = 0$ for $j=0,1,\ldots,m$, $\tilde{a}_{10} = \epsilon$, $\tilde{a}_{i0} = 0$, $i = 2,\ldots,m$, $\tilde{a}_{11} = a_{11} - \epsilon$, $\tilde{a}_{ij} = a_{ij}$ otherwise. It can be regarded as a new coupled network obtained by adding a new node $x_0(t)$ and an edge from the node $x_0(t)$ to node $x_1(t)$ to the original coupled network. Then, the pinning control problem becomes to synchronize the new coupled master-slave network (29). Namely, $\lim_{t\to\infty} |x_i(t) - x_j(t)| = 0$, $i = 0,1,\ldots,m$. Thus, the pinning control problem can be viewed as to synchronize a master-slave network or a connected directed graph with a root node $x_0(t)$ (see [12], [13]).

Moreover, it is clear that all the theoretical results apply to adding more than one controllers. In need, we can pin a coupled network by more than one controllers. Therefore, we can give the following proposition.

*Proposition 2:* Suppose that hypothesis (11) is satisfied. If the corresponding directed graph of the coupled network has a spanning tree with the root set denoted by $root$ and the selective single pinning node $i \in root$, then the coupled network (4) can be synchronized to a chaotic homogenous trajectory. In particular, if the intrinsic network is with an undirected connected graph, then such a single controller can be induced on any node in the network.

This proposition implies that even without much knowledge of the topological structure of the coupled network, we still can pin it by a single controller. Despite that more pinning controllers are sure to improve pinning performance, the study for this simplest case lays down the theoretical basis for further investigation.

## VII. Conclusion

Synchronization is an important research field in sciences and applications. How to pin a coupled network to a specified solution (or an equilibrium point) of the uncoupled system is of great significance. However, in practice, the state variables of some nodes are not observable or measured. Therefore, we have to investigate the possibility of pinning a coupled network by adding controllers to those nodes, which can be measured or controlled. In this paper, we prove rigorously that one can pin a complex network by adding a single linear controller to one node with symmetric or asymmetric coupling matrix. We also discuss how to pin the coupled network with nonlinear coupling. How to find suitable coupling strength adaptively is also studied. Simulations also verify our theoretical results.

### Acknowledgment

The authors thank anonymous reviewers's comments, which helped revise the paper, and Prof. G. Chen in City University of Hong Kong for valuable discussions.

### References

[1] S. Boccaletti, V. Latora, Y. Moreno, M. Chavez, and D. U. Hwang, "Complex networks: Structure and dynamics," *Phys. Rep.*, vol. 424, pp. 175–308, 2006.
[2] M. E. J. Newman, "The structure and function of complex networks," *SIAM Rev.*, vol. 45, pp. 167–256, 2003.
[3] R. Albert and A. L. Barabási, "Statistical mechanics of complex networks," *Rev. Mod. Phys.*, vol. 74, pp. 47–97, 2002.
[4] "Special issue on nonlinear waves, patterns, and spatiotemporal chaos in dynamical arrays," *IEEE Trans. Circuits Syst. I, Fundam. Theory Appl.*, vol. 42, no. 10, pp. 559–828, Oct. 1995.
[5] L. M. Pecora and T. L. Carroll, "Master stability functions for synchronized coupled systems," *Phys. Rev. Lett.*, vol. 80, no. 10, pp. 2109–2112, 1998.
[6] L. M. Pecora, T. Carroll, G. Johnson, D. Mar, and K. S. Fink, "Synchronization stability in coupled oscillator arrays: Solution for arbitrary configuration," *Int. J. Bifurc. Chaos*, vol. 10, no. 2, pp. 273–290, 2000.
[7] X. F. Wang and G. Chen, "Synchronization in scale-free dynamical networks: Robustness and fragility," *IEEE Trans. Circuits Syst. I, Fundam. Theory Appl.*, vol. 49, no. 1, pp. 54–62, Jan. 2002.
[8] M. Barahona and L. M. Pecora, "Synchronization in small-world systems," *Phys. Rev. Lett.*, vol. 89, no. 5, pp. 1–4, 2002.
[9] T. Nishikawa, A. E. Motter, Y. C. Lai, and F. C. Hoppensteadt, "Heterogeneity in oscillator networks: Are smaller world easier to synchronize?," *Phys. Rev. Lett.*, vol. 91, no. 1, pp. 1–4, 2003.
[10] C. W. Wu and L. O. Chua, "Synchronization in an array of linearly coupled dynamical systems," *IEEE Trans. Circuits Syst. I, Fundam. Theory Appl.*, vol. 42, no. 8, 430-447, 1995.
[11] W. L. Lu and T. P. Chen, "Synchronization of coupled connected neural networks with delays," *IEEE Trans. Circuits Syst. I, Reg. Papers*, vol. 51, no. 12, pp. 2491–2503, Dec. 2004.
[12] ——, "New approach to synchronization analysis of linearly coupled ordinary differential systems," *Physica D*, vol. 213, pp. 214–230, 2006.
[13] C. W. Wu, "Synchronization in networks of nonlinear dynamical systems coupled via a directed graph," *Nonlinearity*, vol. 18, pp. 1057–1064, 2005.
[14] C. W. Wu, T. Yang, and L. O. Chua, "On adaptive synchronization and control of nonlinear dynamical systems," *Int. J. Bifurc. Chaos*, vol. 6, no. 3, pp. 455–471, Mar. 1996.
[15] Y. G. Hong, H. S. Qin, and G. R. Chen, "Adaptive synchronization of chaotic systems via state or output feedback control," *Int. J. Bifur. Chaos*, vol. 11, no. 4, pp. 1149–1158, Apr. 2001.
[16] K. Y. Lian, P. Liu, T. S. Chiang, and C. S. Chiu, "Adaptive synchronization design for chaotic systems via a scalar driving signal," *IEEE Trans. Circuits Syst. I*, vol. 49, pp. 17–27, Jan. 2002.
[17] H. G. Z. Qu, "Controlling spatiotemporal chaos in coupled map lattice systems," *Phys. Rev. Lett.*, vol. 72, no. 1, pp. 68–71, 1994.
[18] N. Parekh, S. Parthasarathy, and S. Sinha, "Global and local of spatiotemporal chaos in coupled map lattices," *Phys. Rev. Lett.*, vol. 81, no. 7, pp. 1401–1404, 1998.
[19] R. O. Grigoriev, M. C. Cross, and H. G. Schuster, "Pinning control of spatiotemporal chaos," *Phys. Rev. Lett.*, vol. 79, no. 15, pp. 2795–2798, 1997.
[20] X. F. Wang and G. R. Chen, "Pinning Control of scale-free dynamical-networks," *Phys. A*, vol. 310, pp. 521–531, 2002.
[21] X. Li, X. F. Wang, and G. R. Chen, "Pinning a complex dynamical networks to its equilibrium," *IEEE Transactions on Circuits and Systems-I, Regular Papers*, vol. 51, no. 10, pp. 2074–2087, 2004.
[22] A. Greilich, M. Markus, and E. Goles, "Control of spatiotemporal chaos: Dependence of the minimun pinning distance on the spatial measure entropy," *Eur. Phys. J. D*, vol. 33, pp. 279–283, 2005.

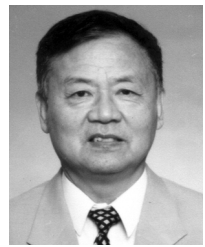

**Tianping Chen** is a Professor in the Department of Mathematics, Fudan University, Shanghai, China. His research interests include harmonic analysis, approximation theory, neural networks, signal processing, dynamical systems and complex networks.

Dr. Chen is a recipient of several important awards, including second prize of 2002' National natural sciences award of China, 1997' outstanding paper award of IEEE Transactions on Neural Networks, 1997' Best paper award of Japanese Neural Network Society.

**Xiwei Liu** is working toward the Ph.D. degree in the School of Mathematical Sciences, Fudan University, Shanghai, China.

His research interests cover nonlinear dynamical systems, complex networks, and neural networks.

**Wenlian Lu** is presently a Postdoctoral Researcher at the Max Planck Institute for Mathematics in the Sciences, Leipzig, Germany, and also a Lecturer in the School of Mathematical Sciences, Fudan University, Shanghai, China. His research interests cover nonlinear dynamical systems, complex networks, and neural networks.